\documentclass{article}

\usepackage{arxiv}

\usepackage[utf8]{inputenc} 
\usepackage[T1]{fontenc}    
\usepackage{hyperref}       
\usepackage{url}            
\usepackage{booktabs}       
\usepackage{amsfonts}       
\usepackage{nicefrac}       
\usepackage{microtype}      
\usepackage{lipsum}		
\usepackage{graphicx}
\usepackage{natbib}
\usepackage{doi}
\usepackage{amsmath, amsthm, amssymb}

\newcommand{\norm}[1]{\left\Vert#1\right\Vert}

\newcommand{\abs}[1]{\left\vert#1\right\vert}
\newcommand{\set}[1]{\left\{#1\right\}}
\newcommand{\seq}[1]{\left<#1\right>}
\newcommand{\bra}[1]{\left(#1\right)}
\newcommand{\sbra}[1]{\left[#1\right]}

\newcommand{\R}{\mathbb{R}}
 \newcommand{\h}{\mathcal{H}}

   \newcommand{\N}{\mathbb{N}}
      \newcommand{\V}{\mathcal{V}}
   \newcommand{\T}{\mathbb{T}}

   \newcommand{\x}{\mathbf{x}}

\renewcommand{\c}{\mathbb{C}}
\newcommand{\n}{\mathbf{\ker}}

\newcommand{\bh}{\mathcal{B(\mathcal{H})}}

\newcommand{\hh}{\mathcal{H}\oplus \mathcal{H}}
\newcommand{\Dst}{\Delta_{s, t}}
\newtheorem{theorem}{Theorem}[section]
\newtheorem{lemma}[theorem]{Lemma}

\newtheorem{corollary}[theorem]{Corollary}
\newtheorem{definition}[theorem]{Definition}

\newtheorem{remark}[theorem]{Remark}

\newcommand\mystyle{\everymath{\displaystyle}}
\mystyle
\title{Elevating Precision in Inequalities for Numerical Radii and Operator Matrices}

\author{ \href{https://orcid.org/0000-0002-3816-5287}{\includegraphics[scale=0.06]{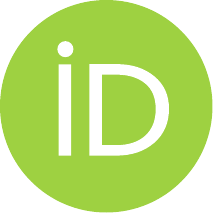}\hspace{1mm}M.H.M.~Rashid}\thanks{Corresponding Author} \\
	Department of Mathematics\&Statistics\\Faculty of Science P.O.Box(7)\\
	Mu'tah University University\\
	Mu'tah-Jordan \\
	\texttt{mrash@mutah.edu.jo}
}

\renewcommand{\shorttitle}{\textit{arXiv} Template}

\hypersetup{
pdftitle={Elevating Precision in Inequalities for Numerical Radii and Operator Matrices},
pdfsubject={q-bio.NC, q-bio.QM},
pdfauthor={M.H.M.Rashid},
pdfkeywords={Cone metric space, Normal cones, Non-normal cones, Fixed point theorem},
}

\begin{document}
\maketitle

\begin{abstract}
	In this paper, we aim to establish a range of numerical radius inequalities. These discoveries will bring us to a recently validated numerical radius inequality and will present numerical radius inequalities that exhibit enhanced precision when compared to those recently established for particular cases. Additionally, we employ the generalized Aluthge transform for operators to deduce a set of inequalities pertaining to the numerical radius. Moreover, we set forth various upper and lower bounds for the numerical radius of $2\times 2$ operator matrices, refining and expanding upon the bounds determined previously.
\end{abstract}

\keywords{Numerical Range\and Numerical Radius\and Aluthge transformation\and Operator matrix}

\section{Introduction}
Let $\h$ be complex Hilbert space and $\bh$ be the $C^*$-algebra   of all bounded linear operator on $\h$.
An operator $A\in\bh$ is said to be {\it positive} if $\seq{Ax,x}\geq 0$ holds for all $x\in\h$. We write $A\geq 0$ if $A$ is positive. \\
\indent  The numerical radius of operator $A$ is formally defined as:
\begin{equation*}
  w(A)=\sup\{\abs{\seq{Ax,x}}:x\in\h,\norm{x}=1\}.
\end{equation*}
It is widely recognized that $w(\cdot)$ serves as a norm on $\bh$, which is tantamount to the conventional operator norm $\norm{\cdot}$. As a matter of fact, for any $A\in \bh$, the following holds:
\begin{equation}\label{N1}
\frac{1}{2}\norm{A}\leq w(A)\leq \norm{A}.
\end{equation}
(see \cite[Theorem 1.3-1]{GR}). Therefore, the standard operator norm and the numerical radius are essentially equivalent. The inequalities in (\ref{N1}) are precise. The initial inequality becomes an equality when $A^2=0$, while the second inequality becomes an equality when $A$ is normal. For fundamental properties of the numerical radius, please refer to \cite{GR}. Notably, Kittaneh has significantly refined the inequalities in (\ref{N1}) in \cite{kit1.5} and \cite{kit2}.
One significant inequality concerning $w(T)$ is the power inequality, which states that for any natural number $n$, the following holds:
\begin{equation}\label{Ref.1}
  w(T^n)\leq (w(T))^n.
\end{equation}
It has been demonstrated in \cite{kit1.5} and \cite{kit2} that if $A\in\bh$, then the subsequent inequalities are valid:
\begin{equation}\label{N2}
  w(A)\leq \frac{1}{2}\norm{|A|+|A^*|}\leq \frac{1}{2}\bra{\norm{A}+\norm{A^2}^{\frac{1}{2}}},
\end{equation}
where $|A|=\sqrt{A^*A}$ is the absolute value of $A$, and
\begin{equation}\label{N3}
 \frac{1}{4}\norm{A^*A+AA^*}\leq w^2(A)\leq \frac{1}{2}\norm{A^*A+AA^*}.
\end{equation}
The representation of the direct sum of two instances of $\h$ is denoted by $\h\oplus\h$. If we consider $A, B, C, D \in \bh$, the matrix operator
$\begin{bmatrix}A &B \\C & D \\ \end{bmatrix}$ can be seen as an operator that operates on $\h\oplus\h$. Its operation is defined as
\begin{equation*}
\begin{bmatrix}A &B \\C & D \\ \end{bmatrix}z=\begin{bmatrix}Ax+By \\Cx+Dy \\\end{bmatrix}
\end{equation*}
 for all
$z=\begin{pmatrix} x \\ y \\ \end{pmatrix}\in \h\oplus \h$. Therefore, the inner product on $\h\oplus\h$ is defined by
\begin{equation*}
  \seq{\begin{bmatrix}A &B \\C & D \\ \end{bmatrix}z,z}=\seq{\begin{bmatrix}Ax+By \\Cx+Dy \\\end{bmatrix},\begin{bmatrix} x \\ y \\ \end{bmatrix}}
  =\seq{Ax+By,x}+\seq{Cx+Dy,y}.
\end{equation*}
For further information, readers are advised to consult \cite{BP1}, \cite{BBP2},\cite{HK}, \cite{Rash1} and \cite{RB}.\\
\indent The paper is organized as follows: Section 2 is dedicated to establishing a variety of numerical radius inequalities. Among the outcomes, we illustrate that for $A,B\in\bh$, where $B=X+iY$ is the Cartesian decomposition of $B$, the following inequality holds
\begin{equation*}
  w^2(AB\pm BA^*)\leq 2\norm{A}^2\bra{w(X^2)+w(Y^2)+\sqrt{\bra{w(X^2)-w(Y^2)}^2+w^2(XY+YX)}}.
\end{equation*}
 Section 3 is dedicated to examining the concept of the generalized Aluthge transform for operators, leading to a series of inequalities concerning the numerical radius. Moreover, we present upper bounds for the numerical radius, particularly tailored for the context of $2\times 2$ operator matrices.
In the final section, section 4, we establish a variety of upper and lower bounds for the numerical radius of $2\times 2$ operator matrices. The derived bounds enhance and generalize the previously known bounds. Specifically, we demonstrate if $A,B,X,Y\in\bh$, then
\begin{equation*}
  w(X^*AY+Y^*BX)\leq 2\norm{X}\norm{Y}w\bra{\begin{bmatrix} 0 &A \\ B& 0 \\\end{bmatrix}}.
\end{equation*}
In particular,
$$w(X^*AY+Y^*AX)\leq 2\norm{X}\norm{Y}w\bra{A}.$$
\section{Numerical radius for sums and products of operator}
In this section, we aim to establish a range of numerical radius inequalities. These findings will guide us to a recently proven numerical radius inequality, and we will introduce numerical radius inequalities that exhibit greater precision compared to those recently established for specific cases.\\
\indent In order to demonstrate our extended numerical radius inequalities, we rely on several widely recognized lemmas. The first lemma, known as the H\"older-McCarthy inequality, is an established outcome derived from the spectral theorem for positive operators and Jensen's inequality (see \cite{kit1}).
\begin{lemma}(H\"older Mc-Carty inequality).\label{Holder}
  Let $A\in\bh$, $A\geq 0$  and let $x\in\h$ be any unit vector. Then  we have
  \begin{enumerate}
    \item [(i)] $\seq{Ax,x}^r\leq \seq{A^rx,x}$ for $r\geq 1$.
    \item [(ii)] $\seq{A^rx,x}\leq \seq{Ax,x}^r$ for $0<r\leq 1$.
  \end{enumerate}
\end{lemma}
The second lemma directly stems from the well-known Jensen's inequality, which addresses the convexity or concavity of specific power functions. It represents a specific instance of Schl\"omilch's inequality for weighted means of non-negative real numbers.
\begin{lemma}\label{J}
  Let $a,b>0$ and $0\leq \alpha\leq 1$. Then
  \begin{equation}
    a^{\alpha}b^{1-\alpha}\leq \alpha a+(1-\alpha)b\leq \left(\alpha a^r+(1-\alpha)b^r\right)^{\frac{1}{r}}\quad\mbox{for}\quad r\geq 1.
  \end{equation}
\end{lemma}
The ensuing lemma results from the convex nature of the function $f(t)=t^r$, where $r\geq 1$ (refer to \cite{Bha, Bohr, AA}).
\begin{lemma}(Bohr's inequality)\label{Logain1}
  Let $a_i,i=1,\cdots,n$ be positive real numbers. Then
  \begin{equation}\label{hoopy1}
    \left(\sum_{i=1}^{n}a_i\right)^r\leq n^{r-1}\sum_{i=1}^{n}a_i^r\quad \mbox{for}\,\,r\geq 1.
  \end{equation}
\end{lemma}
The subsequent lemma holds significant utility in the subsequent discussion, commonly referred to as the generalized mixed Schwartz inequality (refer to \cite{kit1}).
\begin{lemma}\label{Lem2.3}
  Let $T\in\bh$ and $x,y\in\h$ be any vectors.
  \begin{enumerate}
    \item [(i)] If $\alpha,\beta\geq 0$ such that $\alpha+\beta=1$, then $|\seq{Tx,y}|^2\leq \seq{|T|^{2\alpha}x,x}\seq{|T^*|^{2\beta}y,y}$,
    \item [(ii)] If $f,g$ are non-negative continuous functions on $[0,\infty)$  satisfying $f(t)g(t)=t$ ($t\geq 0$),
    $|\seq{Tx,y}|\leq \norm{f(|T|)x}\norm{g(|T^*|)y}$.
  \end{enumerate}
\end{lemma}
\begin{theorem}\label{Rahma1}
  Let $X ,Y \in\bh$ and let $f,g$
   be nonnegative continuous function on $[0,\infty)$  and satisfy the relation $f(t)g(t)=t$ for all $t\geq 0$. Then
  \begin{equation}\label{Boshra1}
    \abs{\seq{Xx,x}}^r+\abs{\seq{Yx,x}}^s\leq \leq \frac{1}{2}\norm{f^{2r}(|X|)+g^{2r}(|X^*|)+f^{2s}(|Y|)+g^{2s}(|Y^*|)}
  \end{equation}
  for all $x\in\h$ and $r,s\geq 1$.
\end{theorem}
\begin{proof} For any $x\in\h$, we have
\begin{eqnarray*}
  &&\abs{\seq{Xx,x}}^r+\abs{\seq{Yx,x}}^s\leq \sbra{\norm{f(|X|x)}\norm{g(|X^*|x)}}^{r}+\sbra{\norm{f(|Y|x)}\norm{g(|Y^*|x)}}^{s}\\
  &&\bra{\mbox{by Lemma \ref{Lem2.3}}}\\
  &&\leq \sbra{\norm{f(|X|x)}^2\norm{g(|X^*|x)}^2}^{\frac{r}{2}}+\sbra{\norm{f(|Y|x)}^2\norm{g(|Y^*|x)}^2}^{\frac{s}{2}}\\
  &&\leq \sbra{\seq{f^2(|X|)x,x}^{\frac{1}{2}}\seq{g^2(|X^*|)x,x}^{\frac{1}{2}}}^{r}+
  \sbra{\seq{f^2(|Y|)x,x}^{\frac{1}{2}}\seq{g^2(|Y^*|)x,x}^{\frac{1}{2}}}^{s}\\
  &&\leq \sbra{\frac{1}{2}\seq{f^2(|X|)x,x}+\frac{1}{2}\seq{g^2(|X^*|)x,x}}^{r}+
  \sbra{\frac{1}{2}\seq{f^2(|Y|)x,x}+\frac{1}{2}\seq{g^2(|Y^*|)x,x}}^{s}\\
  &&\bra{\mbox{by AG-Mean}}\\
  &&\leq \sbra{\frac{1}{2}\seq{f^2(|X|)x,x}^r+\frac{1}{2}\seq{g^2(|X^*|)x,x}^r}+
  \sbra{\frac{1}{2}\seq{f^2(|Y|)x,x}^s+\frac{1}{2}\seq{g^2(|Y^*|)x,x}^s}\\
  &&\bra{\mbox{by Lemma \ref{J}}}\\
  &&\leq \sbra{\frac{1}{2}\seq{f^{2r}(|X|)x,x}+\frac{1}{2}\seq{g^{2r}(|X^*|)x,x}}+
  \sbra{\frac{1}{2}\seq{f^{2s}(|Y|)x,x}+\frac{1}{2}\seq{g^{2s}(|Y^*|)x,x}}\\
  &&\bra{\mbox{by Lemma \ref{Holder}}}\\
  &&=\frac{1}{2}\seq{\bra{f^{2r}(|X|)+g^{2r}(|X^*|)+f^{2s}(|Y|)+g^{2s}(|Y^*|)}x,x}\\
  &&\leq \frac{1}{2}\norm{f^{2r}(|X|)+g^{2r}(|X^*|)+f^{2s}(|Y|)+g^{2s}(|Y^*|)}
\end{eqnarray*}
as desired.
\end{proof}
\begin{corollary} Let $A\in\bh$. Then
\begin{equation*}
  w(A^2)\leq \frac{1}{2}\norm{|A^2|+|A^{2*}|}.
\end{equation*}
\end{corollary}
\begin{proof} In inequality (\ref{Boshra1}), letting $X=Y=A^2$, $s=r=1$ and $f(t)=g(t)=t^{\frac{1}{2}}$ , we obtain
\begin{equation*}
  \abs{\seq{A^2x,x}}\leq \frac{1}{2}\norm{|A^2|+|A^{2*}|}.
\end{equation*}
Taking the supremum over all vectors $x\in\h$ with $\norm{x}=1$, we get
\begin{equation*}
  w(A^2)\leq \frac{1}{2}\norm{|A^2|+|A^{2*}|}.
\end{equation*}
\end{proof}
\begin{corollary} Let $A,B\in\bh$. Then
\begin{equation*}
  w^r(B^*A)\leq \frac{1}{2}\norm{|B^*A|^{r}+|A^*B|^{r}}.
\end{equation*}
\end{corollary}
\begin{proof} The proof follows immediately from Theorem \ref{Rahma1} by setting
$X=Y=B^*A$, $r=s$ and $f(t)=g(t)=\sqrt{t}$.
\end{proof}
\begin{theorem}\label{ALQADRI1} Let $A,B\in\bh$ and let let $f,g$
   be nonnegative continuous function on $[0,\infty)$  and satisfy the relation $f(t)g(t)=t$ for all $t\geq 0$. Then
   \begin{eqnarray}\label{Micro1}
   && \abs{\seq{Ax,x}}^{r}\abs{\seq{Bx,x}}^{s}\leq\\
    && \frac{1}{2}\norm{\alpha f^{\frac{2r}{\alpha}}(|A|)+
     \alpha g^{\frac{2r}{\alpha}}(|A^*|)+(1-\alpha)f^{\frac{2s}{1-\alpha}}(|B|)+(1-\alpha)g^{\frac{2s}{1-\alpha}}(|B^*|)}\nonumber
   \end{eqnarray}
  for all $x\in\h$, $r,s\geq 1$ and $0<\alpha<1$.
\end{theorem}
\begin{proof} For every $x\in\h$, we have
\begin{eqnarray*}
  && \abs{\seq{Ax,x}}^{r}\abs{\seq{Bx,x}}^{s}\leq \sbra{\norm{f(|A|)x}\norm{g(|A^*|x)}}^{r} \sbra{\norm{f(|B|)x}\norm{g(|B^*|x)}}^{s}\\
  &&\bra{\mbox{by Lemma \ref{Lem2.3}}}\\
 &&=\bra{\sbra{\seq{f^2(|A|)x,x}^{\frac{1}{2}}\seq{g^2(|A^*|)x,x}^{\frac{1}{2}}}^{\frac{r}{\alpha}}}^{\alpha}
 +\bra{\sbra{\seq{f^2(|B|)x,x}^{\frac{1}{2}}\seq{g^2(|B^*|)x,x}^{\frac{1}{2}}}^{\frac{s}{1-\alpha}}}^{1-\alpha} \\
 &&\leq \alpha\sbra{\seq{f^2(|A|)x,x}^{\frac{1}{2}}\seq{g^2(|A^*|)x,x}^{\frac{1}{2}}}^{\frac{r}{\alpha}}+
 (1-\alpha)\sbra{\seq{f^2(|B|)x,x}^{\frac{1}{2}}\seq{g^2(|B^*|)x,x}^{\frac{1}{2}}}^{\frac{s}{1-\alpha}}\\
 &&\bra{\mbox{by Lemma \ref{J}}}\\
 &&\leq \frac{\alpha}{2}\sbra{\seq{(f^{\frac{2r}{\alpha}}(|A|)+g^{\frac{2r}{\alpha}}(|A^*|))x,x}}+
 \frac{1-\alpha}{2}\sbra{\seq{(f^{\frac{2s}{1-\alpha}}(|B|)+g^{\frac{2s}{1-\alpha}}(|B^*|))x,x}}\\
 &&\bra{\mbox{by Lemma \ref{J} and Lemma \ref{Holder}}}\\
 &&=\frac{1}{2}\seq{\bra{\alpha f^{\frac{2r}{\alpha}}(|A|)+\alpha g^{\frac{2r}{\alpha}}(|A^*|)
 +(1-\alpha)f^{\frac{2s}{1-\alpha}}(|B|)+(1-\alpha) g^{\frac{2s}{1-\alpha}}(|B^*|)}x,x}\\
 &&\leq \frac{1}{2}\norm{\alpha f^{\frac{2r}{\alpha}}(|A|)+\alpha g^{\frac{2r}{\alpha}}(|A^*|)
 +(1-\alpha)f^{\frac{2s}{1-\alpha}}(|B|)+(1-\alpha) g^{\frac{2s}{1-\alpha}}(|B^*|)}
\end{eqnarray*}
as required.
\end{proof}
\begin{remark}  In Theorem \ref{ALQADRI1}, by setting $r=s$, $A=B$, $\alpha=\frac{1}{2}$, and $f(t)=g(t)=\sqrt{t}$, we arrive at the following result:
\begin{equation*}
  \abs{\seq{Ax,x}}^{2r}\leq \frac{1}{2}\norm{|A|^{2r}+|A^*|^{2r}}.
\end{equation*}
By considering the supremum over all $x\in\h$ with $\norm{x}=1$, we obtain
\begin{equation}\label{F1}
  w^{2r}(A)\leq \frac{1}{2}\norm{|A|^{2r}+|A^*|^{2r}}.
\end{equation}
In inequality (\ref{F1}), if we take $r=1$, we get
\begin{equation*}
  w^{2}(A)\leq \frac{1}{2}\norm{|A|^{2}+|A^*|^{2}}.
\end{equation*}
\end{remark}
\begin{definition} Let $T\in\bh$. Then the number $c(T)$ is defined as follows:
\begin{equation*}
  c(T):=\inf\set{\abs{\seq{Tx,x}}:x\in\h, \norm{x}=1}.
\end{equation*}
\end{definition}
\begin{lemma}\cite{MP} \label{MohdEl}Let $a,b,e\in\h$. Then we have
\begin{eqnarray}\label{FG1}
 &&\abs{\seq{a,b}}\leq \sqrt{\bra{\norm{a}^2-\seq{a,e}^2}}\sqrt{\norm{b}^2-\seq{b,e}^2}+\abs{\seq{a,e}\seq{e,b}} \\
  &&\leq \min\set{\norm{a}\sqrt{\bra{\norm{a}^2-\seq{a,e}^2}},\norm{b}\sqrt{\norm{b}^2-\seq{b,e}^2}}+\abs{\seq{a,e}\seq{e,b}}.\nonumber
\end{eqnarray}
\end{lemma}
\begin{theorem} Let $A,B\in\bh$. Then
\begin{eqnarray*}
  w(B^*A) &\leq& \sqrt{\norm{A}^2+c^2(A)}\sqrt{\norm{B}^2+c^2(B)}+ w(A)w(B) \\
   &\leq&\min\set{\norm{A}\sqrt{\norm{A}^2+c^2(A)},\norm{B}\sqrt{\norm{B}^2+c^2(B)}}+w(A)w(B).
\end{eqnarray*}
In particular, if $A=B$, then
\begin{equation*}
  w(A^2)\leq \norm{A}^2+c^2(A)+w^2(A).
\end{equation*}
\end{theorem}
\begin{proof} The demonstration readily ensues from Lemma \ref{MohdEl} by selecting $a=Ax$, $b=Bx$, and $e=x$, followed by computing the supremum over all vectors $x\in\h$ with $\norm{x}=1$.
\end{proof}
\begin{theorem}\label{Theorem2.13} Let $A,B\in\bh$, and let $B=X+iY$ be the Cartesian decomposition of $B$. Then
\begin{equation*}
  w^2(AB\pm BA^*)\leq 2\norm{A}^2\bra{w(X^2)+w(Y^2)+\sqrt{\bra{w(X^2)-w(Y^2)}^2+w^2(XY+YX)}}.
\end{equation*}
\end{theorem}
\begin{proof} By using the simple identity $Re\bra{AB+BA^*}=ARe(B)+Re(B)A^*$, we have
\begin{eqnarray*}
  \abs{\seq{Re\bra{e^{i\theta}\bra{AB+BA^*}}x,x}}^2 &\leq&\abs{\seq{\bra{ARe(e^{i\theta}B)+Re(e^{i\theta}B)A^*}x,x}}^2\\
  &&\bra{\mbox{since $\abs{\seq{Re(T)x,x}}\leq \abs{\seq{Tx,x}}$}} \\
   &\leq&2\bra{\abs{\seq{ARe(e^{i\theta}B)x,x}}^2+\abs{\seq{Re(e^{i\theta}B)A^*x,x}}^2}\\
   &&\bra{\mbox{by the triangle inequality and the convexity of $f(t)=t^2$}}\\
   &=&  2\bra{\abs{\seq{Re(e^{i\theta}B)x,A^*x}}^2+\abs{\seq{A^*x,Re(e^{i\theta}B)x}}^2}\\
   &\leq& 2\norm{A^*}^2\bra{\norm{Re(e^{i\theta}B)x}^2+\norm{Re(e^{i\theta}B)x}^2}\\
   &=& 4\norm{A}^2\seq{Re(e^{i\theta}B)x,Re(e^{i\theta}B)x}\\
   &=&4\norm{A}^2\seq{\bra{Re(e^{i\theta}B)}^2x,x}.
\end{eqnarray*}
It follows from
\begin{eqnarray*}
  \bra{Re(e^{i\theta}B)}^2 &=& \bra{Re(e^{i\theta}(X+iY))}^2 \\
   &=& \bra{\cos\theta X-\sin\theta Y}^2=\cos^2\theta X^2+\sin^2\theta Y^2-(XY+YX)\sin\theta\cos\theta
\end{eqnarray*}
that
\begin{eqnarray*}
 && \sup_{\theta\in\R}\seq{\bra{Re(e^{i\theta}B)}^2x,x}
 = \sup_{\theta\in\R}\bra{\cos^2\theta \seq{X^2x,x}+\sin^2\theta \seq{Y^2x,x}-\sin\theta\cos\theta\seq{(XY+YX)x,x}}\\
   &\leq&\sup_{\theta\in\R}\bra{\cos^2\theta w(X^2)+\sin^2\theta  w(Y^2)-\sin\theta\cos\theta\seq{(XY+YX)x,x}}\\
   &\leq& \frac{1}{2}\bra{w(X^2)+w(Y^2)+\sqrt{\bra{w(X^2)-w(Y^2)}^2+\bra{\seq{(XY+YX)x,x}}^2}}.
\end{eqnarray*}
Hence,
\begin{eqnarray*}
  &&\sup_{\theta\in\R}\abs{\seq{Re\bra{e^{i\theta}\bra{AB+BA^*}}x,x}}^2  \\
  &&\leq  2\norm{A}^2\bra{w(X^2)+w(Y^2)+\sqrt{\bra{w(X^2)-w(Y^2)}^2+\bra{\seq{(XY+YX)x,x}}^2}}.
\end{eqnarray*}
It is known that if $T\in\bh$, then
\begin{equation*}
  \abs{\seq{Tx,x}}=\sup_{\theta\in\R}Re\bra{e^{i\theta}\seq{Tx,x}}\,\,\mbox{and}\,\, w(T)=\sup_{\theta\in\R}\norm{Re\bra{e^{i\theta}T}}
  =\sup_{\theta\in\R} w\bra{e^{i\theta}T}.
\end{equation*}
Consequently, by taking the supremum over all unit vectors $x\in\h$, we have
\begin{eqnarray}\label{Masafa-1}
  w^2(AB+BA^*) &=&\sup_{\theta\in\R}w\bra{e^{i\theta}\bra{AB+BA^*}} \\
   &\leq& 2\norm{A}^2\bra{w(X^2)+w(Y^2)+\sqrt{\bra{w(X^2)-w(Y^2)}^2+w^2(XY+YX)}}.\nonumber
\end{eqnarray}
Replacing $A$ by $iA$ in the inequality (\ref{Masafa-1}), we get
\begin{equation*}
  w^2(AB-BA^*)\leq 2\norm{A}^2\bra{w(X^2)+w(Y^2)+\sqrt{\bra{w(X^2)-w(Y^2)}^2+w^2(XY+YX)}}.
\end{equation*}
Therefore,
\begin{equation*}
  w^2(AB\pm BA^*)\leq 2\norm{A}^2\bra{w(X^2)+w(Y^2)+\sqrt{\bra{w(X^2)-w(Y^2)}^2+w^2(XY+YX)}}.
\end{equation*}
\end{proof}
Theorem \ref{Theorem2.13} includes a special case as follows.
\begin{corollary} Let $A,B\in\bh$, and let $B=X+iY$ be the Cartesian decomposition of $B$. Then
\begin{enumerate}
  \item [(i)] If $XY+YX=0$, then $w(AB\pm BA^*)\leq 2\norm{A}\max\set{w^{\frac{1}{2}}(X^2),w^{\frac{1}{2}}(Y^2)}$.
  \item [(ii)] If $B$ is self-adjoint, then  $w(AB\pm BA^*)\leq 2\norm{A}w^{\frac{1}{2}}(X^2)$.
  \item [(iii)]  If $B$ is self-adjoint, then $w(AB)\leq \norm{A}w^{\frac{1}{2}}(X^2)$.
\end{enumerate}
\end{corollary}
\begin{proof} (i) The  inequality follows from Theorem \ref{Theorem2.13} and the inequality
\begin{eqnarray*}
  &&w^2(AB\pm BA^*)\leq 2\norm{A}^2\bra{w(X^2)+w(Y^2)+\sqrt{\bra{w(X^2)-w(Y^2)}^2}} \\
  &&=2\norm{A}^2\bra{w(X^2)+w(Y^2)+\abs{w(X^2)-w(Y^2)}}\\
  &&=4\norm{A}^2\max\set{w(X^2),w(Y^2)}.
\end{eqnarray*}
(ii) The second inequality follows from Theorem \ref{Theorem2.13} and the hypothesis $B=X+0i$.\\
(iii) We have
\begin{eqnarray*}
  &&w(AB)=\sup_{\theta\in\R}w\bra{Re\bra{e^{i\theta}AB}} \\
  &&= \frac{1}{2}\sup_{\theta\in\R}w\bra{e^{i\theta}AB+e^{-i\theta}BA^*}\\
  &&\leq \norm{A}w^{\frac{1}{2}}(X^2)\,\,\bra{\mbox{by part (ii)}}.
\end{eqnarray*}
So, the proof is complete.
\end{proof}
\begin{theorem} Let $A,B\in\bh$. Then
\begin{enumerate}
  \item [(i)] $w(AB\pm BA)\leq w^{\frac{1}{2}}(A^*A+AA^*)w^{\frac{1}{2}}(B^*B+BB^*)$.
  \item [(ii)] $w(AB\pm BA)\leq w^{\frac{1}{2}}(A^*A+B^*B)w^{\frac{1}{2}}(AA^*+BB^*)$.
\end{enumerate}
\end{theorem}
\begin{proof} Let $x$ be any unit vector in $\h$. Then
\begin{eqnarray*}
  \abs{\seq{(AB\pm BA)x,x}} &\leq&\abs{\seq{ABx,x}}+\abs{\seq{BAx,x}} \\
   &=& \abs{\seq{Bx,A^*x}}+\abs{\seq{Ax,B^*x}}\\
   &\leq& \norm{Bx}\norm{A^*x}+\norm{Ax}\norm{B^*x}\\
   &\leq& \bra{\norm{A^*x}^2+\norm{Ax}^2}^{\frac{1}{2}}\bra{\norm{B^*x}^2+\norm{Bx}^2}^{\frac{1}{2}}\\
   &&\bra{\mbox{by the Cauchy-Schwarz inequality}}\\
   &=&\abs{\seq{AA^*x,x}+\seq{A^*Ax,x}}^{\frac{1}{2}}\abs{\seq{BB^*x,x}+\seq{B^*Bx,x}}^{\frac{1}{2}}\\
   &=& \abs{\seq{(AA^*+A^*A)x,x}}^{\frac{1}{2}}\abs{\seq{(BB^*+B^*B)x,x}}^{\frac{1}{2}}\\
   &\leq& w^{\frac{1}{2}}(A^*A+AA^*)w^{\frac{1}{2}}(B^*B+BB^*).
\end{eqnarray*}
Hence,
\begin{equation*}
  w(AB\pm BA)=\sup_{\norm{x}=1}\abs{\seq{(AB\pm BA)x,x}} \leq w^{\frac{1}{2}}(A^*A+AA^*)w^{\frac{1}{2}}(B^*B+BB^*).
\end{equation*}
Now, according to the inequality
\begin{eqnarray*}
  &&  \norm{Bx}\norm{A^*x}+\norm{Ax}\norm{B^*x}\leq \bra{\norm{A^*x}^2+\norm{B^*x}^2}^{\frac{1}{2}}\bra{\norm{Ax}^2+\norm{Bx}^2}^{\frac{1}{2}}
 \\
  &&\bra{\mbox{by the Cauchy-Schwarz inequality}}.
\end{eqnarray*}
and a similar argument of the proof of part (i), we get the second inequality.
\end{proof}
For the special case $A=I$, we have the next result.
\begin{corollary} Let $B\in\bh$. Then
\begin{enumerate}
  \item [(i)] $w^2(B)\leq \frac{1}{2}w(B^*B+BB^*)$.
  \item [(ii)] $w^2(B)\leq \frac{1}{4}w(I+B^*B)w(I+BB^*)$.
\end{enumerate}
\end{corollary}
\section{Upper bounds for the numerical radius}
In this section, we utilize the generalized Aluthge transform for operators to deduce a series of inequalities related to the numerical radius. Furthermore, we provide upper bounds for the numerical radius, specifically tailored to the realm of $2\times 2$ operator matrices. Our results broaden the scope of various findings documented in the current literature, as evidenced by references such as \cite{JKP}, \cite{Rash1}, \cite{RB}, and \cite{T}.\\
\indent For any operator $T \in \bh$, it is possible to express it as $T=U|T|$, where $U$ is a partial isometry, and $|T|$ represents the square root of $T^*T$. If the uniqueness of $U$ is determined by the kernel condition $\n(U)=\n(|T|)$, then this particular decomposition is referred to as the "polar decomposition," which is a fundamental result in operator theory. In this paper, we use the notation $T=U|T|$ to denote the polar decomposition, satisfying the kernel condition $\n(U)=\n(|T|)$.\\
\indent For non-negative values of $s$ and $t$ satisfying $s+t=1$, we define the $(s,t)$-Aluthge transform as follows:
\begin{equation*}
\Dst(T)=|T|^sU|T|^t,\quad T\in\bh.
\end{equation*}
\begin{theorem}\label{thm2}
  Let $T\in\bh$ and $h$ be non-negative nondecreasing continuous convex function on $[0,\infty)$
  and let $s,t\geq 0$ such that $s+t=1$. Then
  $$h(w(T))\leq \dfrac{1}{2}\left(h(w(\Dst(T))+\norm{h(|T|)}\right).$$
\end{theorem}
To prove Theorem \ref{thm2}, we prepare the following results:
\begin{lemma}(\cite{AF})\label{lem7}
  Let $A,B,C,D\in \bh$. Then
  $$r(AB+CD)\leq \frac{1}{2}(w(BA)+w(DC))+\frac{1}{2}\sqrt{(w(BA)-w(DC))^2+4\norm{BC}\norm{DA}}.$$
\end{lemma}
\begin{lemma}(\cite{Y1})\label{lem8}
Let $T\in\bh$. Then
\begin{enumerate}
  \item [(i)] $w(T)=\displaystyle{\max_{\theta\in\R}\norm{Re(e^{i\theta}T)}}$.
  \item [(ii)] $w\left(\begin{bmatrix} 0 & T \\0& 0 \\ \end{bmatrix}\right)=\frac{1}{2}\norm{T}$
\end{enumerate}
\end{lemma}
\begin{proof}[Proof of Theorem \ref{thm2}]
  Let $T=U|T|$ be the polar decomposition of $T$. Then for every $\theta\in\R$, we have
  \begin{align}\label{Eq.C1}
    \norm{Re(e^{i\theta}T)} & =r(Re(e^{i\theta}T))=\dfrac{1}{2}r\left(e^{i\theta}T+e^{-i\theta}T^*\right)\nonumber \\
     &=\dfrac{1}{2}r\left(e^{i\theta}U|T|+e^{-i\theta}|T|U^*\right)\nonumber\\
     &=\dfrac{1}{2}r\left(e^{i\theta}U|T|^s|T|^t+e^{-i\theta}|T|^t|T|^sU^*\right).
  \end{align}
  If we put $A=e^{i\theta}U|T|^s, B=|T|^t, C=e^{-i\theta}|T|^t$ and $D=|T|^sU^*$ in Lemma \ref{lem7}, then we obtain
  \begin{align}\label{Eq.C2}
    r\left(e^{i\theta}U|T|^s|T|^t+e^{-i\theta}|T|^t|T|^sU^*\right) & \leq \dfrac{1}{2}\left(w(|T|^sU|T|^t)+w(|T|^tU^*|T|^s\right)\nonumber \\
     &+ \dfrac{1}{2}\sqrt{4\norm{e^{-i\theta}|T|^s|T|^t}\norm{|T|^tU^*e^{i\theta}|T|^s}}\,\,(\text{by Lemma \ref{lem7}})\nonumber\\
     &\leq w(|T|^sU|T|^t)+\sqrt{\norm{|T|}\norm{|T|}}\nonumber\\
     &=w(\Dst(T))+\norm{T}.
  \end{align}
  Using inequalities (\ref{Eq.C1}), (\ref{Eq.C2}) and Lemma \ref{lem8} we have
  $$w(T)=\displaystyle{\max_{\theta\in\R}\norm{Re(e^{i\theta}T)}}\leq \dfrac{1}{2}\left(w(\Dst(T)+\norm{T}\right).$$
  Hence
  \begin{align*}
    h(w(T)) &\leq h\left(\dfrac{1}{2}\left(w(\Dst(T)+\norm{T}\right)\right)\,\,(\text{by the monotonicity of $h$ })\\
     &\leq \dfrac{1}{2}h\left(w(\Dst(T))\right)+\dfrac{1}{2}h(\norm{T})\,\,(\text{by the convexity of $h$ }).
  \end{align*}
  This achieves the proof.
\end{proof}
By making use of Theorem \ref{thm2}, we derive the subsequent outcome.
\begin{corollary}\label{cor3}
  Let $T\in\bh$ and $h$ be non-negative nondecreasing continuous convex function on $[0,\infty)$
  and let $s,t\geq 0$ such that $s+t=1$. Then
  $$h(w(T))\leq \dfrac{1}{4}\norm{h(|T|^{2t})+h(|T|^{2s})}+\dfrac{1}{2}h(w(\Dst(T))).$$
\end{corollary}
\begin{proof}
  First, we note that by the hypotheses we have
  \begin{align}\label{Eq.C3}
    h(|T|)& =h(|T|^s|T|^t)\nonumber \\
     &\leq h\left(\frac{|T|^{2s}+|T|^{2t}}{2}\right)\,\,(\text{by the arithmetic-geometric inequality})\nonumber\\
     &\leq \dfrac{1}{2}\left(h(|T|^{2t})+h(|T|^{2s})\right)\,\,(\text{by the convexity of $h$}).
  \end{align}
  Hence, using Theorem \ref{thm2} and inequality (\ref{Eq.C3}) we obtain
  \begin{align*}
    h(w(T)) &\leq \dfrac{1}{2}\left(h(w(\Dst(T)))+\norm{h(|T|)}\right) \\
     &\leq \dfrac{1}{2}\left(h(w(\Dst(T)))+\dfrac{1}{2}\norm{h(|T|^{2t})+h(|T|^{2s})}\right) \\
     &=\dfrac{1}{2}h(w(\Dst(T)))+\dfrac{1}{4}\norm{h(|T|^{2t})+h(|T|^{2s})}.
  \end{align*}
\end{proof}
\begin{theorem}\label{thm3}
  Let $T,S\in\bh$. Then for all non-negative non-decreasing convex function $h$ on $[0,\infty)$
   and all $s,t\geq 0$ such that $s+t=1$, we have
  \begin{align*}
  h\left( w\left(\begin{bmatrix} 0 & T \\S& 0 \\ \end{bmatrix}\right)\right) & \leq \dfrac{1}{4}\max\{\norm{h(|T|^{2s})+h(|T|^{2t})},\norm{h(|S|^{2t})+h(|S|^{2s})}\} \\
     & +\dfrac{1}{4}\left(h(\norm{|S|^s|T^*|^t})+h(\norm{|T|^s|S^*|^t})\right).
  \end{align*}
\end{theorem}
\begin{proof}
  Let $T=U|T|$ and $S=V|S|$  be the polar decompositions of $T$ and $S$, respectively, and let
  $R=\begin{bmatrix} 0 & T \\S& 0 \\ \end{bmatrix}$. It follows from the polar decomposition of
  $R=\begin{bmatrix} 0 & U \\V& 0 \\ \end{bmatrix}\begin{bmatrix} |S| & 0 \\0& |T| \\ \end{bmatrix}$ that
  \begin{align*}
    \Dst(R) &=|R|^s\begin{bmatrix} 0 & U \\V& 0 \\ \end{bmatrix}|R|^t \\
     &=\begin{bmatrix} |S|^s & 0 \\0& |T|^s \\ \end{bmatrix}\begin{bmatrix} 0 & U \\V& 0 \\ \end{bmatrix}\begin{bmatrix} |S|^t & 0 \\0& |T|^t \\ \end{bmatrix}\\
     &=\begin{bmatrix} 0 & |S|^sU|T|^t \\|T|^sV|S|^t& 0 \\ \end{bmatrix}
  \end{align*}
  Using $|T^*|^2=TT^*=U|T|^2U^*$ and $S^*=SS^*=V|S|^2S^*$ we have $|T|^t=U^*|T^*|^tU$ and $|S|^t=V|S^*|^tV^*$. Therefore
  \begin{align}\label{Eq.C4}
    w(\Dst(R)) &=w\left(\begin{bmatrix} 0 & |S|^sU|T|^t \\|T|^sV|S|^t& 0 \\ \end{bmatrix}\right)\nonumber \\
     &\leq w\left(\begin{bmatrix} 0 & |S|^sU|T|^t \\0& 0 \\ \end{bmatrix}\right)+w\left(\begin{bmatrix} 0 & 0 \\|T|^sV|S|^t& 0 \\ \end{bmatrix}\right)\nonumber\\
     &=w\left(\begin{bmatrix} 0 & |S|^sU|T|^t \\0& 0 \\ \end{bmatrix}\right)+w\left(W^*\begin{bmatrix} 0 & 0 \\|T|^sV|S|^t& 0 \\ \end{bmatrix}W\right),\nonumber\\
     &\,\,(\text{where $W=\begin{bmatrix} 0 & I \\I& 0 \\ \end{bmatrix}$ is unitary})\nonumber\\
     &=w\left(\begin{bmatrix} 0 & |S|^sU|T|^t \\0& 0 \\ \end{bmatrix}\right)+w\left(\begin{bmatrix} 0 & |T|^sV|S|^t \\0& 0 \\ \end{bmatrix}\right)\nonumber\\
     &=\dfrac{1}{2}\norm{|S|^sU|T|^t}+\dfrac{1}{2}\norm{|T|^sV|S|^t}\,\,(\text{by (ii) of Lemma \ref{lem8}})\nonumber\\
     &=\dfrac{1}{2}\norm{|S|^sUU^*|T^*|^tU}+\dfrac{1}{2}\norm{|T|^sVV^*|S^*|^tV}\nonumber\\
     &\leq \dfrac{1}{2}\norm{|S|^t|T^*|^t}+\dfrac{1}{2}\norm{|T|^s|S^*|^t}
  \end{align}
  Applying Corollary \ref{cor3} and inequality (\ref{Eq.C4}), we have
  \begin{align*}
             h(w(R)) &\leq \dfrac{1}{4}\norm{h(|R|^{2t})+h(|R|^{2s})}+\dfrac{1}{2}h(w(\Dst(R))) \\
             &\leq \dfrac{1}{4}\max\{\norm{h(|T|^{2t})+h(|T|^{2s})},\norm{h(|S|^{2t})+h(|S|^{2s})}\}\\
             &+\dfrac{1}{2}\left(\dfrac{1}{2}h((\norm{|S|^s|T^*|^t}+\norm{|T|^s|S^*|^t})\right)) \\
             &\leq \dfrac{1}{4}\max\{\norm{h(|T|^{2t})+h(|T|^{2s})},\norm{h(|S|^{2t})+h(|S|^{2s})}\}\\
             &+\dfrac{1}{4}h(\norm{|S|^s|T^*|^t})+h(\norm{|T|^s|S^*|^t})\,\,(\text{by the convexity of $h$}).
 \end{align*}
 This completes the proof.
\end{proof}
\begin{corollary}\label{cor2}
  Let $T,S\in \bh$ and let $s,t\geq 0$ such that $s+t=1$. Then
  $$w^p\left(\begin{bmatrix} 0 & T \\S& 0 \\ \end{bmatrix}\right)\leq \dfrac{1}{2}\max\{\norm{T}^p,\norm{S}^p\}
  +\dfrac{1}{4}\left(\norm{|S|^s|T^*|^t}^{p}+\norm{|T|^s|S^*|^s}^{p}\right),$$
  where $p\geq 1$.
\end{corollary}
\begin{proof}
  By using Theorem \ref{thm2} and inequality (\ref{Eq.C4}), we have
  \begin{align*}
    w^p\left(\begin{bmatrix} 0 & T \\S& 0 \\ \end{bmatrix}\right) & \leq \dfrac{1}{2}\norm{\begin{bmatrix} 0 & T \\S& 0 \\ \end{bmatrix}}+\dfrac{1}{2}w^p(\Dst(R)) \\
     &=\dfrac{1}{2}\max\{\norm{T},\norm{S}\}+\dfrac{1}{2}\left(\dfrac{1}{2}(\norm{|S|^s|T^*|^t}+\norm{|T|^s|S^*|^t})\right)^{p}\\
     &\leq \dfrac{1}{2}\max\{\norm{T},\norm{S}\}+\dfrac{1}{4}(\norm{|S|^s|T^*|^t}^{p}+\norm{|T|^s|S^*|^t}^{p})
  \end{align*}
  this completes the proof.
\end{proof}

\begin{corollary}\label{cor4}
  Let $T,S\in\bh$. Then, for all $s,t\geq 0$ such that $s+t=1$ and $p\geq 1$, we have
  \begin{align*}
    w^{\frac{p}{2}}(TS) &\leq \dfrac{1}{4}\max\left\{\norm{|T|^{2sp}+|T|^{2tp}},\norm{|S|^{2sp}+|S|^{2tp}}\right\} \\
     &+\dfrac{1}{4}\{\norm{|T|^s|S^*|^t}^p+\norm{|S|^s|T^*|^t}^p\}.
  \end{align*}
\end{corollary}
\begin{proof}
  Applying the power inequality of the numerical radius, we get
  \begin{align*}
     w^{\frac{p}{2}}(TS) &\leq \max\{ w^{\frac{p}{2}}(TS), w^{\frac{p}{2}}(ST)\} \\
     &= w^{\frac{p}{2}}\left(\begin{bmatrix} TS & 0 \\0& ST \\ \end{bmatrix}\right)=w^{\frac{p}{2}}\left(\begin{bmatrix} 0 & T \\S& 0 \\ \end{bmatrix}^2\right)\\
     &\leq w^p\left(\begin{bmatrix} 0 & T \\S& 0 \\ \end{bmatrix}\right)\\
     &\leq \dfrac{1}{4}\max\left\{\norm{|T|^{2sp}+|T|^{2tp}},\norm{|S|^{2sp}+|S|^{2tp}}\right\} \\
     &+\dfrac{1}{4}\{\norm{|T|^s|S^*|^t}^p+\norm{|S|^s|T^*|^t}^p\}\,\,(\text{by Theorem \ref{thm3}}),
  \end{align*}
  so the proof is complete.
\end{proof}
\begin{corollary}
  Let $T,S\in\bh$ be positive operators. Then, for all $s,t\geq 0$ such that $s+t=1$ and $p\geq 1$, we have
  $$\norm{T^{\frac{1}{2}}S^{\frac{1}{2}}}^p\leq \dfrac{1}{4}\max\{\norm{T^{sp}+T^{tp}},\norm{S^{sp}+S^{tp}}\}
  +\dfrac{1}{4}\{\norm{T^sS^t}^p+\norm{S^sT^t}^p\}.$$
\end{corollary}
\begin{proof}
  Since the spectral radius of any operator is dominated by its numerical radius, we have
  $\sqrt{r(TS)}\leq \sqrt{w(TS)}$. Applying a commutativity property of the spectral radius, we have
  \begin{align*}
    r^{\frac{p}{2}}(TS) & =r^{\frac{p}{2}}(T^{\frac{1}{2}}T^{\frac{1}{2}}S^{\frac{1}{2}}S^{\frac{1}{2}}) \\
    &r^{\frac{p}{2}}(T^{\frac{1}{2}}S^{\frac{1}{2}}S^{\frac{1}{2}}T^{\frac{1}{2}})\\
    &=r^{\frac{p}{2}}((T^{\frac{1}{2}}S^{\frac{1}{2}})(T^{\frac{1}{2}}S^{\frac{1}{2}})^*)\\
    &=\norm{T^{\frac{1}{2}}S^{\frac{1}{2}}}^{p}.
      \end{align*}
      So, the result follows now by Corollary \ref{cor4}.
\end{proof}
\begin{corollary}\label{cor4.9}
  Let $T,S\in\bh$. Then for all non-negative non-decreasing convex function $h$ on $[0,\infty)$
   and all $s,t\geq 0$ such that $s+t=1$, we have
   \begin{align*}
     h(\norm{T+S}) &\leq \dfrac{1}{4}\max\{\norm{h(2|T|^{2s})+h(2|T|^{2t})}, \norm{h(2|S^*|^{2s})+h(2|S^*|^{2s})}\}\\
      &+\dfrac{1}{4}\{h(2\norm{|T^*|^s|S|^t})+h(2\norm{|S^*|^s|T^*|^t})\}.
   \end{align*}
   In particular, if $T$ and $S$ are normal, then
   $$ h(\norm{T+S})\leq \dfrac{1}{4}\max\{h(2\norm{|T|}),h(2\norm{|S|})\}+\dfrac{1}{4}h(2\norm{TS}^{\frac{1}{2}}). $$
\end{corollary}
\begin{proof}
  Let $R=\begin{bmatrix} 0 & T \\S& 0 \\ \end{bmatrix}$. Applying Lemma \ref{lem8} and Corollary \ref{cor3}, we have
  \begin{align*}
    h(\norm{T+S^*}) &=h(\norm{R+R^*}) \\
     &\leq \displaystyle{\max_{\theta\in\R} h\left(\norm{2Re(e^{i\theta}R)}\right)}\\
     &=h(2w(R))\\
     &\leq \dfrac{1}{4}\max\{\norm{h(2|T|^{2s})+h(2|T|^{2t})},\norm{h(2|S|^{2s})+h(2|S|^{2t})}\}\\
     &+\dfrac{1}{4}\{h(2\norm{|T|^s|S^*|^t})+h(2\norm{|S|^s|T^*|^t})\}\,\,(\text{by Theorem \ref{thm3}}).
  \end{align*}
  Now, the desired result follows by replacing $S$ by $S^*$. For the particular case $t=s=\frac{1}{2}$. Since
  $T$ and $S$ are normal, we have $|T|=|T^*|$ and $|S|=|S^*|$ and
  \begin{align*}
    h(\norm{|T|^{\frac{1}{2}}|S|^{\frac{1}{2}}}) &=h(\sqrt{r(|T||S|)}) \\
     & \leq h(\norm{|T||S|}^{\frac{1}{2}})\,\,(\text{by the convexity of $h$})\\
     &=h(\norm{U^*TS^*V}^{\frac{1}{2}})\\
     &=h(\norm{TS^*}^{\frac{1}{2}}),
  \end{align*}
  where $T=U|T|$ and $S=V|S|$ are the polar decomposition of $T$ and $S$. This completes the proof.
\end{proof}
\begin{corollary}
  Let $T,S\in\bh$ and $p\geq 1$. Then
  \begin{align*}
    \norm{T+S}^{p} & \leq \dfrac{1}{2^{2-r}}\max\{\norm{|T|^{2tp}+|T|^{2sp}},\norm{|S^*|^{2tp}+|S^*|^{2sp}}\} \\
     & +\dfrac{1}{2^{2-r}}\left\{\norm{|T|^t|S|^s}^{p}+\norm{|S^*|^t|T^*|^s}^p\right\}.
  \end{align*}
  In particular, if $T$ and $S$ are normal, then
  $$\norm{T+S}^p\leq \dfrac{1}{2^{1-r}}\max\{\norm{T}^p,\norm{S}^p\}+\dfrac{1}{2^{1-r}}\norm{TS}^{\frac{p}{2}}.$$
\end{corollary}
\begin{proof}
  The proof follows from Corollary \ref{cor4.9} for the convex function $h(x)=x^p$ $(p\geq 1)$.
\end{proof}
We end this section with the following result.
\begin{theorem}\label{Theorem 3.11} Let $A,B\in\bh$ and let $\nu\in [0,1]$. Then
\begin{equation}\label{Ineq.3E1}
  w^{2r}(B^*A)\leq \norm{\nu|A|^{\frac{2r}{\nu}}+(1-\nu)|B|^{\frac{2r}{1-\nu}}}
\end{equation}
for every $r\geq 1$.
\end{theorem}
\begin{proof} For every unit vector $x\in\h$, we have
\begin{eqnarray*}
  \abs{\seq{B^*Ax,x}}^{2r}&=&\abs{\seq{Ax,Bx}}^{2r}\leq \norm{Ax}^{2r}\norm{Bx}^{2r}\,\,\bra{\mbox{by Cauchy-Schwarz inequality}}\\
  &\leq&\seq{|A|^{2r}x,x}\seq{|B|^{2r}x,x} \,\,\bra{\mbox{by Lemma \ref{Holder}}}\\
  &\leq&\seq{|A|^{\frac{2r}{\nu}}x,x}^{\nu}\seq{|B|^{\frac{2r}{1-\nu}}x,x}^{1-\nu} \,\,\bra{\mbox{by Lemma \ref{Holder}}}\\
  &\leq& \nu \seq{|A|^{\frac{2r}{\nu}}x,x}+(1-\nu)\seq{|B|^{\frac{2r}{1-\nu}}x,x}\,\,\bra{\mbox{by Young's Inequality}}\\
  &=&\seq{\bra{\nu |A|^{\frac{2r}{\nu}}+(1-\nu)|B|^{\frac{2r}{1-\nu}} }x,x}\\
  &\leq& \norm{\nu |A|^{\frac{2r}{\nu}}+(1-\nu)|B|^{\frac{2r}{1-\nu}}}.
\end{eqnarray*}
Taking the supremum over all unit vector $x\in\h$, we obtain
\begin{equation*}
  w^{2r}(B^*A)\leq \norm{\nu|A|^{\frac{2r}{\nu}}+(1-\nu)|B|^{\frac{2r}{1-\nu}}}
\end{equation*}
for all $r\geq 1$ and $\nu\in[0,1]$. This completes the proof.
\end{proof}
\section{numerical radius of $2 \times 2$
operator matrices}
In this segment, various upper and lower limits for the numerical radius of $2\times 2$ operator matrices are established, improving and extending the previously identified bounds.
\begin{lemma}\cite{HK}\label{lemma5.1} Let $T,S\in\bh$.  Then the following assertions hold:
\begin{enumerate}
  \item [(a)] $w\bra{\begin{bmatrix}
                       T &0 \\
                       0 & S \\
                     \end{bmatrix}
  }=\max\set{w(T),w(S)}$.
  \item [(b)] $w\bra{\begin{bmatrix}
                       0 &T \\
                       S & 0 \\
                     \end{bmatrix}
  }=w\bra{\begin{bmatrix}
                       0 &S \\
                       T & 0 \\
                     \end{bmatrix}}$.
  \item [(c)] $w\bra{\begin{bmatrix}
                       0 &T \\
                       S & 0 \\
                     \end{bmatrix}}=w\bra{\begin{bmatrix}
                       0 &T \\
                      e^{i\theta} S & 0 \\
                     \end{bmatrix}}$ for all $\theta\in\R$.
\item [(d)] $w\bra{\begin{bmatrix}
                       T &S \\
                       S & T \\
                     \end{bmatrix}}=\max\set{w(T-S),w(T+S)}$. In particular, $w\bra{\begin{bmatrix}
                       0 &S \\
                       S & 0 \\
                     \end{bmatrix}}=w(S)$.
\end{enumerate}
\end{lemma}
Subsequently, we will establish an upper limit for the numerical radius of $2\times 2$ operator matrices.
\begin{theorem}\label{theorem5.2} Let $A,B,C,D\in\bh$. Then
\begin{equation*}
  w^{r}\bra{\begin{bmatrix}
                       A &B \\
                       C& D \\
                     \end{bmatrix}}\leq 4^{r-1}w\bra{\begin{bmatrix}
                       w^{r}(A) &\norm{B}^{r} \\
                       \norm{C}^{r} & w(D)^{r} \\
                     \end{bmatrix}}
\end{equation*}
for all $r\geq 1$.
\end{theorem}
\begin{proof} Suppose $z=(x,y)\in\h\oplus\h$ with $\norm{z}=1$, where $\norm{x}^2+\norm{y}^2=1$. Then
\begin{eqnarray*}
  &&\abs{\seq{\begin{bmatrix}
                       A &B \\
                       C& D \\
                     \end{bmatrix}z,z}}^{2r}=\abs{\seq{\begin{bmatrix}
                       A &B \\
                       C& D \\
                     \end{bmatrix}\begin{pmatrix}
                                    x \\
                                    y \\
                                  \end{pmatrix},\begin{pmatrix}
                                    x \\
                                    y \\
                                  \end{pmatrix}}}^{2r}\\
   &=&\abs{\seq{\begin{pmatrix}
                                   A x+By \\
                                    Cx+Dy \\
                                  \end{pmatrix},\begin{pmatrix}
                                    x \\
                                    y \\
                                  \end{pmatrix}}}^{r} =\abs{\seq{Ax,x}+\seq{By,x}+\seq{Cx,y}+\seq{Dy,y}}^{r}\\
   &\leq& 4^{r-1}\bra{\abs{\seq{Ax,x}}^{r}+\abs{\seq{Dy,y}}^{r}+\abs{\seq{By,x}}^{r}+\abs{\seq{Cx,y}}^{r}}\\
   &\leq& 4^{r-1}\bra{w^{r}(A)\norm{x}^{2r}+w^{r}(D)\norm{y}^{2r}+\norm{B}^{r}\norm{x}^{r}\norm{y}^{r}+\norm{C}^{r}\norm{x}^{r}\norm{y}^{r}}\\
   &=&4^{r-1}\bra{\seq{\begin{bmatrix}
                       w^{r}(A) &\norm{B}^{r} \\
                       \norm{C}^{r}& w(D)^{r} \\
                     \end{bmatrix}\widetilde{z},\widetilde{z}}}\quad\bra{\mbox{where $\widetilde{z}=\bra{\norm{x}^{r},\norm{y}^{r}}\in\c^2$}}.
\end{eqnarray*}
Thus
\begin{equation*}
  w^{r}\bra{\begin{bmatrix}
                       A &B \\
                       C& D \\
                     \end{bmatrix}}=\sup_{\norm{z}=1}\abs{\seq{\begin{bmatrix}
                       A &B \\
                       C& D \\
                     \end{bmatrix}z,z}}^{r}\leq 4^{r-1}w\bra{\begin{bmatrix}
                       w^{r}(A) &\norm{B}^{r} \\
                       \norm{C}^{r} & w(D)^{r} \\
                     \end{bmatrix}}.
\end{equation*}
\end{proof}
Observe that $w(A)\leq \norm{A}$, $w(D)\leq\norm{D}$, and $w([a_{ij }]) \leq w([b_{ij} ])$, for any $0\leq a_{ij }\leq b_{ij }$. Consequently, the subsequent corollary follows directly from Theorem \ref{theorem5.2}.
\begin{corollary}\label{corollary5.3} Let $A,B,C,D\in\bh$. Then
\begin{equation*}
  w^r\bra{\begin{bmatrix}
                       A &B \\
                       C& D \\
                     \end{bmatrix}}\leq 4^{r-1} \begin{bmatrix}
                       \norm{A}^r &\norm{B}^r \\
                       \norm{C}^r & \norm{D}^r \\
                     \end{bmatrix}.
\end{equation*}
\end{corollary}
It is worth noting that Theorem \ref{theorem5.2} provides a more improved bound compared to the one in Corollary \ref{corollary5.3}. To demonstrate the next outcome, we require the subsequent lemma.
\begin{lemma}\label{lemma5.4}\cite{Halmos} Suppose $B=[b_{ij}]$ is an $n\times n$ matrix where $b_{ij} \geq 0$ for all $i,j= 1, 2,\cdots,n$. Then $w(B)=r\bra{\frac{b_{ij}+b_{ji}}{2}}$, where $r(\cdot)$ represents the spectral radius.
\end{lemma}
By employing Theorem \ref{theorem5.2} and Corollary \ref{corollary5.3} with $r=1$, and leveraging Lemma \ref{lemma5.4}, we derive the subsequent two corollaries.
\begin{corollary}\label{corollary5.5}Let $A,B,C,D\in\bh$. Then
\begin{equation*}
  w\bra{\begin{bmatrix}
                       A &B \\
                       C& D \\
                     \end{bmatrix}}\leq r\bra{[c_{ij}]}\leq \frac{1}{2}\bra{w(A)+w(D)+\sqrt{\bra{w(A)-w(D)}^2+\bra{\norm{B}+\norm{C}}^2}},
\end{equation*}
where $c_{11}=w(A),c_{12}=c_{21}=\frac{\norm{B}+\norm{C}}{2}$ and $c_{22}=w(D)$.
\end{corollary}
\begin{corollary}\label{corollary5.6}Let $A,B,C,D\in\bh$. Then
\begin{equation*}
  w\bra{\begin{bmatrix}
                       A &B \\
                       C& D \\
                     \end{bmatrix}}\leq r\bra{[c_{ij}]}\leq \frac{1}{2}\bra{\norm{A}+\norm{D}+\sqrt{\bra{\norm{A}-\norm{D}}^2+\bra{\norm{B}+\norm{C}}^2}},
\end{equation*}
where $c_{11}=\norm{A},c_{12}=c_{21}=\frac{\norm{B}+\norm{C}}{2}$ and $c_{22}=\norm{D}$.
\end{corollary}
Subsequently, by employing the power inequality (\ref{Ref.1}), we establish a lower limit for the numerical radius of $2\times 2$ operator matrices.
\begin{theorem}\label{Theorem5.7} Let $A,B\in\bh$. Then
\begin{equation*}
  w\bra{\begin{bmatrix}
                       0 &A \\
                       B& 0 \\
                     \end{bmatrix}}\geq \sqrt[2n]{\max\set{(AB)^n,(BA)^n}}.
\end{equation*}
\end{theorem}
\begin{proof} Let $T=\begin{bmatrix}
                       0 &A \\
                       B& 0 \\
                     \end{bmatrix}$. Then $T^{2n}=\begin{bmatrix}
                       (AB)^{n} &0 \\
                       0& (BA)^n \\
                     \end{bmatrix}$ for all $n\in\N$. Using Lemma \ref{lemma5.1} and inequality (\ref{Ref.1}), we acquire:
\begin{equation*}
  \max\set{w((AB)^n),w((BA)^n)}=w(T^{2n})\leq w^{2n(T)},
\end{equation*}
Thus, the outcome is as follows.
\end{proof}
In the subsequent steps, we demonstrate the subsequent lower and upper limits.
\begin{theorem}\label{Theorem5.8} Let $A,B\in\bh$. Then
\begin{equation*}
  \frac{1}{2}\max\set{w(A-B),w(A+B)}\leq w\bra{\begin{bmatrix}
                       0 &A \\
                       B& 0 \\
                     \end{bmatrix}}\leq \frac{1}{2}\bra{w(A-B)+w(A+B)}.
\end{equation*}
\end{theorem}
\begin{proof} From Lemma \ref{lemma5.1}(e), we derive the following:
\begin{eqnarray}\label{Ineq.P1}
  &&w(A+B)=w\bra{\begin{bmatrix}
                       0 &A+B \\
                      A+ B& 0 \\
                     \end{bmatrix}}=w\bra{\begin{bmatrix}
                       0 &A \\
                       B& 0 \\
                     \end{bmatrix}+\begin{bmatrix}
                       0 &B \\
                       A& 0 \\
                     \end{bmatrix}}\nonumber\\
   &\leq&w\bra{\begin{bmatrix}
                       0 &A \\
                       B& 0 \\
                     \end{bmatrix}}+w\bra{\begin{bmatrix}
                       0 &B \\
                       A& 0 \\
                     \end{bmatrix}} = 2w\bra{\begin{bmatrix}
                       0 &A \\
                       B& 0 \\
                     \end{bmatrix}}.
\end{eqnarray}
By replacing $B$ with $-B$, we obtain
\begin{equation}\label{Ineq.P2}
  w(A-B)\leq 2w\bra{\begin{bmatrix}
                       0 &A \\
                       -B& 0 \\
                     \end{bmatrix}}=2w\bra{\begin{bmatrix}
                       0 &A \\
                       B& 0 \\
                     \end{bmatrix}}.
\end{equation}
Consequently, the initial inequality follows from (\ref{Ineq.P1}) and (\ref{Ineq.P2}). To prove the second inequality, we consider a unitary operator
$U=\frac{1}{\sqrt{2}}\begin{bmatrix}
                       I &-I \\
                       I& I \\
                     \end{bmatrix}$. his yields
\begin{eqnarray*}
&&w\bra{\begin{bmatrix} 0 &A \\ B& 0 \\\end{bmatrix}}=w\bra{U^*\begin{bmatrix} 0 &A \\ B& 0 \\\end{bmatrix}U} \\
&&=\frac{1}{2}w\bra{\begin{bmatrix} A+B&A-B \\ -(A-B)& -(A+B) \\\end{bmatrix}} \\
&&\leq \frac{1}{2}w\bra{\begin{bmatrix} A+B &0 \\ 0& -(A+B) \\\end{bmatrix}} +\frac{1}{2}w\bra{\begin{bmatrix} 0 &A-B \\ -(A-B)& 0 \\\end{bmatrix}}\\
&&=\frac{1}{2}\bra{w(A+B)+w(A-B)}.
\end{eqnarray*}
This concludes the proof.
\end{proof}
Applying Theorem \ref{Theorem5.8}, we establish the subsequent inequalities.
\begin{corollary}\label{corollary5.9} Let $T\in\bh$  and let $T=A+iB$ be the cartesian decomposition of $T$. Then
\begin{equation*}
  \frac{1}{2}w(T)\leq w\bra{\begin{bmatrix} 0 &A \\ e^{i\theta}B& 0 \\\end{bmatrix}}\leq w(T).
\end{equation*}
\end{corollary}
\begin{proof} By substituting $B$ with $iB$ in Theorem \ref{Theorem5.8}, and subsequently utilizing Lemma \ref{lemma5.1}, we obtain the following:
\begin{equation*}
  \frac{1}{2}\max\set{w(A+iB),w(A-iB)}\leq w\bra{\begin{bmatrix} 0 &A \\ e^{i\theta}B& 0 \\\end{bmatrix}}\leq \frac{1}{2}\bra{w(A+iB)+w(A-iB)}.
\end{equation*}
This infers that
\begin{equation*}
  \frac{1}{2}\max\set{w(T),w(T^*)}\leq w\bra{\begin{bmatrix} 0 &A \\ e^{i\theta}B& 0 \\\end{bmatrix}}\leq \frac{1}{2}\bra{w(T)+w(T^*)}.
\end{equation*}
However, it is known that $w(T)=w(T^*)$. This concludes the proof.
\end{proof}
Considering the formula for $w\bra{\begin{bmatrix} A &B \\ B& A \\\end{bmatrix}}$ derived in Lemma \ref{lemma5.1}, it is intuitive to seek a comparable expression for $w\bra{\begin{bmatrix} A &B \\ -B& -A \\\end{bmatrix}}$. To establish this, we must first demonstrate the following lemma.
\begin{lemma}\label{lemma5.10} Let $A,B,C,D\in\bh$. Then
\begin{equation*}
  w\bra{\begin{bmatrix} A &B \\ C& D \\\end{bmatrix}}\geq w\bra{\begin{bmatrix} A &0 \\ 0& D \\\end{bmatrix}}\quad
  \mbox{and}\quad w\bra{\begin{bmatrix} A &B \\ C& D \\\end{bmatrix}}\geq w\bra{\begin{bmatrix} 0 &B \\ C& 0 \\\end{bmatrix}}.
\end{equation*}
\end{lemma}
\begin{proof} Let $z=(x,0)\in\h\oplus\h$ with $\norm{z}=1$, that is, $\norm{x}=1$. Consequently, we have:
\begin{equation*}
  \abs{\seq{\begin{bmatrix} A &B \\ C& D \\\end{bmatrix}z,z}}=\abs{\seq{(Ax,Cx),(x,0)}}=\abs{\seq{Ax,x}}.
\end{equation*}
Taking the supremum over $\norm{z}=1$, we arrive at:
\begin{equation*}
  w\bra{\begin{bmatrix} A &B \\ C& D \\\end{bmatrix}}=\sup_{\norm{x}=1}\abs{\seq{Ax,x}}.
\end{equation*}
Consequently, this implies:
\begin{equation}\label{P3}
  w\bra{\begin{bmatrix} A &B \\ C& D \\\end{bmatrix}}\geq w(A).
\end{equation}
Likewise,
\begin{equation}\label{P4}
  w\bra{\begin{bmatrix} A &B \\ C& D \\\end{bmatrix}}\geq w(D).
\end{equation}
Thus, the first inequality is a direct result of (\ref{P3}) and (\ref{P4}), combined with Lemma \ref{lemma5.1} (a). To prove the second inequality, we express:
\begin{equation*}
\begin{bmatrix} 0 &B \\ C& 0 \\\end{bmatrix}=\frac{1}{2}\begin{bmatrix} A &B \\ C& D \\\end{bmatrix}+\frac{1}{2}\begin{bmatrix} -A &B \\ C& -D \\\end{bmatrix}.
\end{equation*}
Consequently
\begin{equation*}
w\bra{\begin{bmatrix} 0 &B \\ C& 0 \\\end{bmatrix}}=\frac{1}{2}w\bra{\begin{bmatrix} A &B \\ C& D \\\end{bmatrix}}
+w\bra{\frac{1}{2}\begin{bmatrix} -A &B \\ C& -D \\\end{bmatrix}}.
\end{equation*}
By considering the unitary operator $U=\begin{bmatrix} 0 &-I \\ I& 0 \\\end{bmatrix}$, we have
$U^*\begin{bmatrix} -A &B \\ C& -D \\\end{bmatrix}U=\begin{bmatrix} -D &-C \\ -B& -A \\\end{bmatrix}$, and utilizing the property
$w(U^*XU)=w(X)$, we deduce that:
\begin{equation}\label{P5}
  w\bra{\begin{bmatrix} 0 &B \\ C& 0 \\\end{bmatrix}}\leq \frac{1}{2}w\bra{\begin{bmatrix} A &B \\ C& D \\\end{bmatrix}}
  +w\bra{\begin{bmatrix} -D &-C \\ -B& -A \\\end{bmatrix}}.
\end{equation}
Again, considering the unitary operator $U=\begin{bmatrix} 0 &I \\ I& 0 \\\end{bmatrix}$, we have
$U^*\begin{bmatrix} -D &-C \\ -B& -A \\\end{bmatrix}U=\begin{bmatrix} -A &-B \\ -C& -D \\\end{bmatrix}$
and $w\bra{\begin{bmatrix} -D &-C \\ -B& -A \\\end{bmatrix}}=w\bra{\begin{bmatrix} -A &-B \\ -C& -D \\\end{bmatrix}}
=w\bra{\begin{bmatrix} A &B \\ C& D \\\end{bmatrix}}$.
By employing this argument, the second desired inequality follows from inequality(\ref{P5}).
\end{proof}
\begin{theorem}\label{theorem5.11}  Let $A,B\in\bh$. Then
\begin{equation*}
  \max\set{w(A),w(B)}\leq w\bra{\begin{bmatrix} A &B \\ -B& -A \\\end{bmatrix}}\leq w(A)+w(B).
\end{equation*}
\end{theorem}
\begin{proof} The first inequality is a consequence of Lemma \ref{lemma5.10} in conjunction with Lemma \ref{lemma5.1}. As for the other part,
 \begin{equation*}
   w\bra{\begin{bmatrix} A &B \\ -B& -A \\\end{bmatrix}}\leq w\bra{\begin{bmatrix} A &0 \\ 0& -A \\\end{bmatrix}}
   +w\bra{\begin{bmatrix} 0 &B \\ -B& 0 \\\end{bmatrix}}=w(A)+w(B).
 \end{equation*}
\end{proof}
Especially, by setting $B=A$ in Theorem \ref{theorem5.11}, we obtain the subsequent inequality.
\begin{corollary}\label{corollary5.12} Let $A\in\bh$. Then
  \begin{equation*}
  w(A)\leq w\bra{\begin{bmatrix} A &A \\ -A& -A \\\end{bmatrix}}\leq 2w(A).
\end{equation*}
\end{corollary}
Subsequently, we establish lower and upper bounds for the numerical radius of $2 \times 2$ operator matrices.
\begin{theorem}Let $A,B,C,D\in\bh$. Then
\begin{eqnarray*}
  &&\max\set{w(A),w(D),\frac{1}{2}w(B+C),\frac{1}{2}w(B-C)}\leq w\bra{\begin{bmatrix} A &B \\ C& D \\\end{bmatrix}} \\
 && \leq \max\set{w(A),w(D)}+\frac{1}{2}\bra{w(B+C)+w(B-C)}.
\end{eqnarray*}
\end{theorem}
\begin{proof} By virtue of Lemma \ref{lemma5.10},
\begin{eqnarray*}
  &&w\bra{\begin{bmatrix} A &B \\ C& D \\\end{bmatrix}} \geq \max\set{w\bra{\begin{bmatrix} A &0 \\ 0& D \\\end{bmatrix}},
  w\bra{\begin{bmatrix} 0 &B \\ C& 0 \\\end{bmatrix}}} \\
  &&=\max\set{w(A),w(D),w\bra{\begin{bmatrix} 0 &B \\ C& 0 \\\end{bmatrix}}} \quad\bra{\mbox{by Lemma \ref{lemma5.1}(a)}}\\
  &&\geq \max\set{w(A),w(D), \frac{1}{2}w(B+C),\frac{1}{2}w(B-C)}\quad\bra{\mbox{by Theorem \ref{Theorem5.8}}}.
\end{eqnarray*}
Once more, we deduce from Lemma \ref{lemma5.1}(a) and Theorem \ref{Theorem5.8} that
\begin{eqnarray*}
  &&w\bra{\begin{bmatrix} A &B \\ C& D \\\end{bmatrix}}\leq w\bra{\begin{bmatrix} A &0 \\ 0& D \\\end{bmatrix}}+
  w\bra{\begin{bmatrix} 0 &B \\ C& 0 \\\end{bmatrix}}\\
 &&\leq \max\set{w(A),w(D)}+\frac{1}{2}\bra{w(B+C)+w(B-C)}.
\end{eqnarray*}
Thus, the proof is concluded.
\end{proof}
Now, utilizing the outcomes derived earlier and employing the identity
\begin{equation}\label{P11}
\max\set{a,b}=\frac{a+b}{2}+\frac{\abs{a-b}}{2}\quad\mbox{for all $a,b\geq 0$},
\end{equation}
we establish the ensuing inequality:
\begin{theorem} Let $A,B\in\bh$. Then
\begin{equation*}
  w\bra{\begin{bmatrix} 0 &A \\ B& 0 \\\end{bmatrix}}+\frac{\abs{w(A+B)-w(A-B)}}{2}\leq w(A)+w(B).
\end{equation*}
\end{theorem}
\begin{proof} Utilizing Theorem \ref{Theorem5.8} and the identity (\ref{P11}), we derive
\begin{eqnarray*}
  && w\bra{\begin{bmatrix} 0 &A \\ B& 0 \\\end{bmatrix}}\leq \frac{w(A+B)+w(A-B)}{2} \\
  &&=\max\set{w(A+B),w(A-B)}-\frac{\abs{w(A+B)-w(A-B)}}{2}\\
  &&\leq w(A)+w(B)-\frac{\abs{w(A+B)-w(A-B)}}{2}.
\end{eqnarray*}
  This establishes the desired inequality.
\end{proof}
The subsequent theorem is stated as follows:
\begin{theorem} Let $A,B,X,Y\in\bh$. Then
\begin{equation*}
  w(X^*AY+Y^*BX)\leq 2\norm{X}\norm{Y}w\bra{\begin{bmatrix} 0 &A \\ B& 0 \\\end{bmatrix}}.
\end{equation*}
In particular,
$$w(X^*AY+Y^*AX)\leq 2\norm{X}\norm{Y}w\bra{A}.$$
\end{theorem}
\begin{proof} Consider $x, y \in \h $ as non-zero vectors, and let $z=\frac{1}{\sqrt{\norm{x}^2+\norm{y}^2}}(x,y)$.
Hence, $z$ is a unit vector in $\h\oplus\h$, leading us to:
\begin{equation*}
  w\bra{\begin{bmatrix} 0 &A \\ B& 0 \\\end{bmatrix}}\geq \abs{\seq{\begin{bmatrix} 0 &A \\ B& 0 \\\end{bmatrix}z,z}}\\
    =\frac{\abs{\seq{Ay,x}+\seq{Bx,y}}}{\norm{x}^2+\norm{y}^2}.
\end{equation*}
 Consequently,
 \begin{equation*}
   \bra{\norm{x}^2+\norm{y}^2}w\bra{\begin{bmatrix} 0 &A \\ B& 0 \\\end{bmatrix}}\geq
   \abs{\seq{Ay,x}+\seq{Bx,y}}.
 \end{equation*}
 Now, upon substituting $x$ and $y$ with $Xx$ and $Yx$, respectively, we obtain
 \begin{equation*}
   \abs{\seq{AYx,Xx}+\seq{BXx,Yx}}\leq \bra{\norm{Xx}^2+\norm{Yx}^2}w\bra{\begin{bmatrix} 0 &A \\ B& 0 \\\end{bmatrix}}.
 \end{equation*}
 Hence
 \begin{equation*}
   \abs{\seq{AYx,Xx}+\seq{BXx,Yx}}\leq w\bra{\begin{bmatrix} 0 &A \\ B& 0 \\\end{bmatrix}}\bra{\norm{X}^2+\norm{Y}^2}\norm{x}^2
 \end{equation*}
 Taking the supremum over $\norm{x}=1$, we obtain
 \begin{equation}\label{P22}
   w(X^*AY+Y^*BX)\leq \bra{\norm{X}^2+\norm{Y}^2}w\bra{\begin{bmatrix} 0 &A \\ B& 0 \\\end{bmatrix}}.
 \end{equation}
 Now, the intended inequality ensues from (\ref{P22}) by substituting $X$ and $Y$ with $tX$ and $\frac{1}{t}Y$ respectively, where $t=\sqrt{\frac{\norm{Y}}{\norm{X}}}$. Notably, by setting $A=B$, we attain the desired second inequality.
\end{proof}
Next, we need the following lemma, known as Buzano’s extension of Schwarz
inequality ( see \cite{Buzano}).
\begin{lemma}\label{Lemma4.16} Let $a,b,e\in\h$ with $\norm{h}=1$. Then
\begin{equation}\label{Ineq.4.1}
  \abs{\seq{a,e}\seq{e,b}}\leq \frac{1}{2}\bra{\norm{a}\norm{b}+\abs{\seq{a,b}}}.
\end{equation}
\end{lemma}
Applying Lemma \ref{J}, we obtain immediately
\begin{lemma}\label{Lemma4.17} Let $a,b,e\in\h$ with $\norm{h}=1$. Then
\begin{equation}\label{Ineq.4.1}
  \abs{\seq{a,e}\seq{e,b}}^{r}\leq \frac{1}{2}\bra{\norm{a}^{r}\norm{b}^{r}+\abs{\seq{a,b}}^{r}}.
\end{equation}
for every $r\geq$.
\end{lemma}
\begin{theorem}Let $B,C\in\bh$. Then
\begin{eqnarray}
   w^{r}\bra{\begin{bmatrix} 0 &A \\ B& 0 \\\end{bmatrix}}&\leq&\frac{1}{4}\max\set{\norm{|C|^{2r}+|B^*|^{2r}},\norm{|B|^{2r}+|C^*|^{2r}}}\nonumber\\
  &+&\frac{1}{2}\max\set{w^r(|B^*||C|),w^r(|C^*||B|)}.
\end{eqnarray}
for every $r\geq 1$
\end{theorem}
\begin{proof} Let $\x$ be a unit vector in $\hh$ and let $\T=\begin{bmatrix} 0 &B \\ C& 0 \\\end{bmatrix}$. Then
\begin{eqnarray*}
  \abs{\seq{\T\x,\x}}^{r} &\leq&\frac{1}{2}\bra{\norm{|\T|\x}^{r}\norm{|\T^*|\x}^{r}+\abs{\seq{|\T|\x,|\T^*|\x}}^{r}}\bra{\mbox{by Lemma \ref{Lemma4.17}}}\\
   &\leq&\frac{1}{4}\bra{\norm{|\T|\x}^{2r}+\norm{|\T^*|\x}^{2r}}+\frac{1}{2}\abs{\seq{|\T|\x,|\T^*|\x}}^{r}\bra{\mbox{by A-G-mean inequality}}\\
   &\leq&\frac{1}{4}\bra{\seq{|\T|^{2r}\x,\x}+\seq{|\T^*|^{2r}\x,\x}}+\frac{1}{2}\abs{\seq{|\T|\x,|\T^*|\x}}^{r}\bra{\mbox{by Lemma \ref{Holder}}} \\
   &\leq&\frac{1}{4}\seq{\begin{bmatrix} |C|^{2r}+|B^*|^{2r} &0 \\ 0& |B|^{2r}+|C^*|^{2r} \\\end{bmatrix}\x,\x}+\frac{1}{2}\abs{\seq{|\T|\x,|\T^*|\x}}^{r}\\
   &\leq&\frac{1}{4}\max\set{\norm{|C|^{2r}+|B^*|^{2r}},\norm{|B|^{2r}+|C^*|^{2r}}}+\frac{1}{2}\max\set{w^r(|B^*||C|),w^r(|C^*||B|)}.
\end{eqnarray*}
By taking the supremum over all vector $\x\in\hh$ with $\norm{\x}=1$, we get
\begin{equation*}
  w^r\bra{\begin{bmatrix} 0 &B \\ C& 0 \\\end{bmatrix}}\leq \frac{1}{4}\max\set{\norm{|C|^{2r}+|B^*|^{2r}},\norm{|B|^{2r}+|C^*|^{2r}}}+\frac{1}{2}\max\set{w^r(|B^*||C|),w^r(|C^*||B|)}.
\end{equation*}
\end{proof}
Using the above lemma we prove the following theorem.
\begin{theorem}\label{Theorem4.18} Let $B,C\in\bh$ and let $\nu\in [0,1]$. Then
  \begin{eqnarray}
    w^{2r}\bra{\begin{bmatrix} 0 &A \\ B& 0 \\\end{bmatrix}}&\leq&\frac{1}{2}\max\set{\norm{\nu |C|^{\frac{2r}{\nu}}+(1-\nu)|B^*|^{\frac{2r}{1-\nu}},\norm{\nu |B|^{\frac{2r}{\nu}}+(1-\nu)|C^*|^{\frac{2r}{1-\nu}}}}}\nonumber \\
   &+& \frac{1}{2}\max\set{w^{2r}(|B^*||C|),w^{2r}(|C^*||B|)}.
  \end{eqnarray}
  for every $r\geq 1$.
\end{theorem}
\begin{proof} Let $\x$ be a unit vector in $\hh$ and let $\T=\begin{bmatrix} 0 &B \\ C& 0 \\\end{bmatrix}$. Then
\begin{eqnarray*}
  \abs{\seq{\T\x,\x}}^{2r}&\leq& \frac{1}{2}\bra{\norm{|\T|\x}^{2r}\norm{|\T^*|\x}^{2r}+|\seq{|\T|\x,|\T^*|\x}|^{2r}}\bra{\mbox{by Lemma \ref{Lemma4.17}}}\\
   &\leq&\frac{1}{2}\bra{\seq{|\T|^{\frac{2r}{\nu}}\x,\x}^{\nu}\seq{|\T^*|^{\frac{2r}{\nu}}\x,\x}^{1-\nu}+|\seq{|\T|\x,|\T^*|\x}|^{2r}}\\
   &\leq&\frac{1}{2}\bra{\seq{\bra{\nu |\T|^{\frac{2r}{\nu}}+(1-\nu)|\T^*|^{\frac{2r}{\nu}}}\x,\x}}+\frac{1}{2}|\seq{|\T|\x,|\T^*|\x}|^{2r}\\
   &&\bra{\mbox{by Young's Inequality}}
\end{eqnarray*}
and so
\begin{eqnarray*}
  \abs{\seq{\T\x,\x}}^{2r} &\leq&\frac{1}{2}\seq{\begin{bmatrix} \nu |C|^{\frac{2r}{\nu}}+(1-\nu)|B^*|^{\frac{2r}{1-\nu}} &0 \\ 0&
   \nu |B|^{\frac{2r}{\nu}}+(1-\nu)|C^*|^{\frac{2r}{1-\nu}}\\\end{bmatrix}}\\
   &+&\frac{1}{2}\abs{\seq{\begin{bmatrix} |B^*||C| &0 \\ 0& |C^*||B| \\\end{bmatrix}}}^{2r} \\
 &\leq&\frac{1}{2}\max\set{\norm{\nu |C|^{\frac{2r}{\nu}}+(1-\nu)|B^*|^{\frac{2r}{1-\nu}},\norm{\nu |B|^{\frac{2r}{\nu}}+(1-\nu)|C^*|^{\frac{2r}{1-\nu}}}}}\\
 &+&\frac{1}{2}\max\set{w^{2r}(|B^*||C|),w^{2r}(|C^*||B|)}.
\end{eqnarray*}
By taking the supremum over all unit vectors $\x\in\hh$, we get
\begin{eqnarray*}
   w^{2r}\bra{\begin{bmatrix} 0 &A \\ B& 0 \\\end{bmatrix}} &\leq&\frac{1}{2}\max\set{\norm{\nu |C|^{\frac{2r}{\nu}}+(1-\nu)|B^*|^{\frac{2r}{1-\nu}},\norm{\nu |B|^{\frac{2r}{\nu}}+(1-\nu)|C^*|^{\frac{2r}{1-\nu}}}}} \\
   &+& \frac{1}{2}\max\set{w^{2r}(|B^*||C|),w^{2r}(|C^*||B|)}.
\end{eqnarray*}
\end{proof}
\begin{remark} (i) If we set $r=1$ and $\nu=\frac{1}{2}$ in Theorem \ref{Theorem4.18}, we obtain
\begin{eqnarray*}
 w^{2}\bra{\begin{bmatrix} 0 &A \\ B& 0 \\\end{bmatrix}} &\leq&\frac{1}{4}\max\set{\norm{ |C|^{4}+|B^*|^{4}},\norm{ |B|^{4}+|C^*|^{4}}} \\
   &+& \frac{1}{2}\max\set{w^2(|B^*||C|),w^2(|C^*||B|)}
\end{eqnarray*}
which obtain by Bhunia and Paul \cite[Theorem 2.10]{BP3}. Consequently, Theorem \ref{Theorem4.18} is an improvement
and generalization of \cite[Theorem 2.10]{BP3}.\\
(ii) In particular, by considering $B=C$ and using Lemma \ref{lemma5.1}, we obtain
\begin{equation}\label{Ineq.4T1}
  w^{2r}(B)\leq \frac{1}{2}\norm{\nu |B|^{\frac{2r}{\nu}}+(1-\nu)|B^*|^{\frac{2r}{1-\nu}}}+\frac{1}{2}w^{2r}(|B^*||B|)
\end{equation}
for every $r\geq 1$ and $\nu\in [0,1]$.\\
(iii) By using Theorem \ref{Theorem 3.11} and the inequality (\ref{Ineq.4T1}), we get
\begin{eqnarray*}
  w^{2r}(B)&\leq& \frac{1}{2}\norm{\nu |B|^{\frac{2r}{\nu}}+(1-\nu)|B^*|^{\frac{2r}{1-\nu}}}+\frac{1}{2}w^{2r}(|B^*||B|) \\
   &\leq& \frac{1}{2}\norm{\nu |B|^{\frac{r}{\nu}}+(1-\nu)|B^*|^{\frac{r}{1-\nu}}}+\frac{1}{2}\norm{\nu |B|^{\frac{2r}{\nu}}+(1-\nu)|B^*|^{\frac{2r}{1-\nu}}}\\
   &\leq& \norm{\nu |B|^{\frac{2r}{\nu}}+(1-\nu)|B^*|^{\frac{2r}{1-\nu}}}.
\end{eqnarray*}
\end{remark}
{\bf Author Contributions:}   The author have read and agreed to the published version of the manuscript.\\
{\bf Funding:} No funding is applicable\\
{\bf Institutional Review Board Statement:} Not applicable.\\
{\bf Informed Consent Statement:} Not applicable.
{\bf Data Availability Statement:} Not applicable.\\
{\bf Conflicts of Interest:} The authors declare no conflict of interest.

\bibliographystyle{unsrtnat}
\bibliography{references}  






\end{document}